# Survival exponents for some Gaussian processes

## G. Molchan

Institute of Earthquake Prediction Theory and Mathematical Geophysics,
Russian Academy of Science, Profsoyuznaya 84/32, Moscow, Russia

*Abstract.* The problem is a power-law asymptotics of the probability that a self-similar process does not exceed a fixed level during long time. The exponent in such asymptotics is estimated for some Gaussian processes, including the fractional Brownian motion (FBM) in $(-T_-,T)$, $T \geq T_- \gg 1$ and the integrated FBM in $(0,T)$, $T \gg 1$.

## 1. The problem

Let $x(t), x(0) = 0$ be a real-valued stochastic process with the following asymptotics:

$$P(x(t) < 1, t \in \Delta_T) = T^{-\theta_x + o(1)}, \quad T \to \infty, \quad (1)$$

where $\theta_x$ is the so-called *survival* exponent of $x(t)$. Below we focus on estimating $\theta_x$ for some self-similar Gaussian processes in extended intervals $\Delta_T = (0,T)$ and $(-T_-,T)$, $T \geq T_- \gg 1$. Usually the estimation of the survival exponents is based on Slepian's lemma. To do this we need reference processes with explicit or almost explicit values of $\theta$. Unfortunately, the list of such processes is very short. This includes the fractional Brownian motion (*FBM*), $w_H(t)$, of order $0 < H < 1$ both with one- and multidimensional time. According to Molchan (1999)

$$\theta_{w_H} = 1 - H \text{ for } \Delta_T = (0,T)) \quad \text{and} \quad \theta_{w_H} = d \text{ for } \Delta_T = (-T,T)^d. \quad (2)$$

Another important example is the integrated Brownian motion $I(t) = \int_0^t w(s)ds$ with the exponent

$$\theta_I = 1/4, \quad \Delta_T = (0,T) \text{ (Sinai, 1992)} \quad (3)$$

The nature of this result is best understood in terms of a series of generalizations where the integrand is random walk with discrete or continuous time (see, e.g., Isozaki and Watanabe, 1994; Isozaki and Kotani, 2000; Simon,2007; Vysotsky,2010; Aurzada and Dereich, 2011; Dembo and Gao, 2011). The extension of (3) to include the case of the integrated fractional Brownian motion, $I_H(t) = \int_0^t w_H(s)ds$, remains an important, but as yet unsolved problem.

Below we consider the survival exponents for the following Gaussian processes:
$I_H(t), t \in (0,T)$; $\chi_H(t) = sign(t) w_H(t)$, $t \in (-T,T)$; *FBM* in $\Delta_T = (-T^\alpha, T)$, $0 \leq \alpha \leq 1$; the Laplace transform of white noise with $\Delta_T = (0,T)$, and the fractional Slepian's stationary process whose correlation function is $B_{S_H}(t) = (1 - |t|^{2H})_+$, $0 < H \leq 1/2$.

Our approach to estimation of $\theta_x$ is more or less traditional. Namely, any self-similar process $x(t)$ in $\Delta_T = (0,T)$ generates a *dual stationary process* $\tilde{x}(s) = e^{-hs} x(e^s)$, $s < \ln T := \tilde{T}$, where $h$ is the self-similarity index of $x(t)$. For a large class of Gaussian processes, relation (1) induces the dual asymptotics

$$P(\tilde{x}(s) \leq 0, 0 < s < \tilde{T}) = \exp(-\tilde{\theta}_x \tilde{T}(1+o(1))), \quad \tilde{T} \to \infty \qquad (4)$$

with the same exponent $\tilde{\theta}_x = \theta_x$ (Molchan, 1999, 2008). More generally, the dual exponent is defined by the asymptotics

$$P(x(t) \leq 0, t \in \Delta_T \setminus (-1,1)) = \exp(-\tilde{\theta}_x \tilde{T}(1+o(1))).$$

To formulate the simplest condition of the exponent equality, we define one more exponent $\breve{\theta}_x$ by means of the asymptotics

$$P(|t_T^*| \leq 1) = T^{-\breve{\theta}_x + o(1)},$$

where $t_T^*$ is the position of the maximum of $x(t)$ in $\Delta_T$, i.e., $x(t_T^*) = \sup(x(t), t \in \Delta_T)$.

**Lemma 1.** Let $x(t), x(0) = 0$ be a self-similar continuous Gaussian process in $\Delta_T = (-T_-, T)$, $T_- \leq T$, and $(H_x(\Delta_T), \|\ \|_x)$ be the reproducing kernel Hilbert space associated with $x(t)$. Suppose there exists such an element $\varphi_T$ of $H_x(\Delta_T)$ that $\varphi_T(t) \geq 1, |t| > 1$ and $\|\varphi_T\|_x^2 = o(\ln T)$. Then $\theta_x, \tilde{\theta}_x$, and $\breve{\theta}_x$ can exist simultaneously only; moreover, the exponents are equal to each other.

The equality $\theta = \tilde{\theta}$ reduces the original problem to the estimation of $\tilde{\theta}$. Non-negativity of the correlation function of $\tilde{x}(s)$ guarantees the existence of the exponent $\tilde{\theta}$ (Li&Shao, 2004). In turn, the inequality of two correlation functions, $B_1(s) \leq B_2(s)$, $B_i(0) = 1$, implies, by Slepian's lemma, the inverse inequality for the corresponding exponents: $\tilde{\theta}_1 \geq \tilde{\theta}_2$.

An essentially different approach is required to find the explicit value of $\theta$ for FBM in $\Delta_T = (-T^\alpha, T)$ and to estimate $\tilde{\theta}$ in (4) for the fractional Slepian process with a small H parameter.

## 2. Examples

### 2.1. Integrated fractional Brownian motion

Consider the process

$$I_H(t) = \int_0^t w_H(s) ds \ ,$$

where $w_H(t)$ is the fractional Brownian motion, i.e., a Gaussian random process with the stationary increments: $E|w_H(t) - w_H(s)|^2 = |t-s|^{2H}$, $w_H(0) = 0$. Molchan and Khokhlov (2003, 2004) analyzed theoretically and numerically the exponent $\theta_{I_H}$ in the general case of $H$ and formulated the following

*Hypothesis*: $\quad \theta_{I_H} = H(1-H)$ for $\Delta_T = (0, T))$ and $\quad \theta_{I_H} = 1-H$ for $\Delta_T = (-T, T)$.

The unexpected symmetry $\theta_{I_H} = \theta_{I_{1-H}}$ of the exponents for $\Delta_T = (0, T)$ caused some doubt as to the numerical results. To support the hypothesis, Molchan(2008) derived the following estimates:

$$\rho H(1-H) \leq \theta_{I_H}^+ \leq \theta_{I_H}^{-/+} \leq (1-H) \ , \qquad (5)$$

where $\rho$ is a small constant and $(+)$ and $(-/+)$ are indicators of the intervals $\Delta_T = (0,T)$ and $\Delta_T = (-T,T)$, respectively. Note that, in the case of $H < 1/2$ and $\Delta_T = (-T,T)$, it is unknown whether the exponent exists. In such cases we have to operate with upper $\overline{\theta}$ and lower $\underline{\theta}$ exponents. Therefore, $\theta_{I_H}^{-/+}$ in (5) for $H < 1/2$ is any number from the interval $(\underline{\theta}, \overline{\theta})$. The relation (5) can be improved as follows:

**Proposition 1.** For the intervals $\Delta_T = (0,T)$

(a) $\theta_{I_H} \geq \theta_{I_{1-H}}$, $0 < H \leq 0.5$,

(b) $0.5(H \wedge \overline{H}) \leq \theta_{I_H} \leq \overline{H}$, $\overline{H} = 1 - H$,

(c) $\theta_{I_H} \leq \sqrt{(1 - (H \wedge \overline{H})^2)/12}$.

**Proof.** The identity of dual exponents for $I_H(t)$ follows from (Molchan and Khokhlov, 2004); the dual survival exponent exists, because the dual correlation function

$$\widetilde{B}_{I_H}(s) = (2+4H)^{-1}[(2+2H)(e^{Hs} + e^{-Hs}) - e^{(1+H)s} - e^{-(1+H)s} + (e^{s/2} - e^{-s/2})^{2H+2}], \quad (6)$$

is positive. The inequality (**a**) is a consequence of the relation

$$\widetilde{B}_{I_H}(t) \leq \widetilde{B}_{I_{1-H}}(t), \quad 0 < H \leq 1/2. \quad (7)$$

To prove (**b, c**), we use the correlation function of the process $\widetilde{I}_{1/2}(ps)$, i.e.,

$$\widetilde{B}_{I_{1/2}}(ps) = 1/2(3\exp(-p|s|/2) - \exp(-3p|s|/2)), \quad (8)$$

and the respective exponent $\widetilde{\theta} = p/4$ (see (3)). The relation

$$\widetilde{B}_{I_H}(t) \leq \widetilde{B}_{I_{1/2}}(pt), \quad H \geq 1/2, \quad p = 2(1-H) \quad (9)$$

implies $\theta_{I_H} \geq (1-H)/2$ for $H \geq 1/2$. Using (**a**) in addition, we come to the lower bound in (**b**) because $\theta_{I_H} \geq \theta_{I_{1-H}} \geq H/2$ for $H \leq 1/2$.

Similarly, the relation

$$\widetilde{B}_{I_H}(t) \geq \widetilde{B}_{I_{1/2}}(pt), \quad H \leq 1/2, \quad p = 2\sqrt{(1-H^2)/3} \quad (10)$$

implies (**c**) for all $H$. A test of the pure analytical facts (7, 9, and 10) is given in Appendix.

**Remark 1.** The proposition 1a follows from the more informative relation

$$P(\widetilde{I}_H(s) \leq 0, s \in (0, \widetilde{T})) \leq P(\widetilde{I}_{1-H}(s) \leq 0, s \in (0, \widetilde{T})). \quad (11)$$

This inequality is important for understanding the numerical result by Molchan and Khokhlov(2003) represented in the form of the empirical estimates of $\widetilde{\theta}_{I_H}$ in Figure 1. We can see that the empirical estimates show small but one-sided deviations from the hypothetical curve $\theta = H(1-H)$ before and after

H=1/2. The signs of these deviations are consistent with (11), while the amplitudes are compatible with the model

$$P(\widetilde{I}_H(s) \leq 0, s \in (0, \widetilde{T})) \approx C\widetilde{T}^{\alpha(H)} \exp(-H(1-H)\widetilde{T}), \quad \widetilde{T} >> 1, \quad \text{sgn } \alpha(H) = sign(H - 0.5), \quad (12)$$

and $\alpha(H) = H - 0.5$ (see more in Molchan and Khokhlov, 2003).

### 2.2 Laplace transform of white noise.

Consider the process $L(t) = t\int_0^\infty e^{-tu} dw(u)$, where $w(u)$ is Brownian motion. The dual stationary process $\widetilde{L}(s)$ has the correlation function $\widetilde{B}_L(s) = 1/\cosh(s/2)$. Using (8) as a majorant of $\widetilde{B}_L(s)$, we improve the lower bound of $\widetilde{\theta}_L$ as follows:

**Proposition 2.** $3^{-1/2} \leq 4\widetilde{\theta}_L \leq 1$.

*Proof*. The exponent equality for the dual processes $L$ and $\widetilde{L}$ follows from Lemma 1 with $\varphi_T(t) = t(1+\varepsilon_T)/(t+\varepsilon_T)$, where $\varepsilon_T = 1/\sqrt{\ln T}$. For indeed, $\varphi_T(t) = EL(t)\eta$, where $\eta = (1+\varepsilon_T^{-1})L(\varepsilon_T)$. By definition of the Hilbert space $H_x(\Delta_T)$, we have the desired estimate:

$$\|\varphi_T\|_L^2 = E\eta^2 = \varepsilon_T^{-1}(\varepsilon_T + 1)^2/2 = O(\sqrt{\ln T}).$$

By (3) and Slepian's lemma, the relation

$$\widetilde{B}_{I_{1/2}}(t) \leq \widetilde{B}_L(pt), \quad p \leq 1 \tag{13}$$

has as a consequence the estimate $4p\widetilde{\theta}_L \leq 1$. The opposite inequality

$$\widetilde{B}_{I_{1/2}}(t) \geq \widetilde{B}_L(pt), \quad p^2 \geq 3 \tag{14}$$

implies $4p\widetilde{\theta}_L \geq 1$. The test of (13, $p = 1$) and (14, $p = 2$) is very simple and yields the Li and Shao (2004) estimates: $0.5 < 4\widetilde{\theta}_L < 1$. The Appendix contains a proof of (13, 14) for all interesting values of $p: 1, 2,$ and $\sqrt{3}$.

**Remark 2.** The *dual survival* exponent of $L(t)$ is of interest as a parameter of the following asymptotic relation

$$P(\sum_0^{2n} \xi_i x^i \neq 0, x \in R^1) = (2n)^{-4\widetilde{\theta}_L + o(1)}, \quad n \to \infty \tag{15}$$

for random polynomials with the standard Gaussian independent coefficients (Dembo et al., 2002). A continuous analogue of the polynomial on any of four intervals $0 < \pm x^{\pm 1} \leq 1$ is the Laplace transform of white noise that partially explains the appearance of $\widetilde{\theta}_L$ in the asymptotic relation (15). Simulations suggest $4\widetilde{\theta}_L = 0.76 \pm 0.03$ (Dembo et al, 2002) and $4\widetilde{\theta}_L \approx 0.75$ (Newman and Loinaz, 2001).

### 2.3 Fractional Slepian's process.

We reserve this term for a Gaussian stationary process $S_H(t)$ with correlation function

$$B_{S_H}(t) = (1-|t|^{2H})_+, \quad 0 < H \leq 1/2, \tag{16}$$

because $S_{1/2}(t)$ is known as the Slepian process and $S_H(t) - S_H(0), 0 < t \leq 1$, is equal in distribution to the fractional Brownian motion on the interval $(0,1)$. By the Polya criterion, the fractional Slepian process exists because $B_{S_H}(t)$ is a non-increasing and convex function on the semi-axis $t \geq 0$. The fact of the correlation function being non-negative guarantees the existence of the exponent $\widetilde{\theta}_{S_H}$ in (4). $S_H(t)$ can be useful as a reference process in estimation of the survival exponents. Therefore it is important to have accurate estimates of the exponent for $S_H(t)$. The case of small $H$ is the most interesting because it describes a transition of $S_H(t)$ to white noise. Our estimates of $\widetilde{\theta}_{S_H}$ are based on two lemmas, where we use the following notation

$$\widetilde{\theta}(f, \Delta) = -|\Delta|^{-1} \log P(x(t) \leq f(t), t \in \Delta). \tag{17}$$

**Lemma 2.** (Li and Shao, 2004). Let $x(t)$ be a centered Gaussian stationary process with a finite non-negative correlation function, i.e., $B_x(t) \geq 0$ and $B_x(t) = 0$ for $|t| \geq T_0$. Then the limit

$$\widetilde{\theta}(a) = \lim_{T \to \infty} \widetilde{\theta}(a, (0, T))$$

exists for every $a \in R^1$. Moreover,

$$(1 + 1/k)^{-1} \widetilde{\theta}(a, k\Delta_0) \leq \widetilde{\theta}(a) \leq \widetilde{\theta}(a, k\Delta_0), \qquad \Delta_0 = (0, T_0). \tag{18}$$

**Remark 3.** Lemma 1 was derived by Li and Shao (2004) for the Slepian process, $S_{1/2}(t)$, but the proof remains valid for the general case. There is an explicit but very complicated formula for $\widetilde{\theta}_{S_H}(0, \Delta)$ with $H = 1/2$ (Shepp, 1971). In case of $\Delta = (0, 2)$, the Shepp result reduces to

$$P(S_{1/2}(t) \leq 0, t \in (0,2)) = 1/6 - (2 + \sqrt{3})/(8\pi),$$

and gives $1.336 < \widetilde{\theta}_{S_{1/2}} < 2.004$.

**Lemma 3.** (Aurzada&Dereich, 2011). Let $x(t)$ be a centered Gaussian process in an interval $\Delta$ with a correlation function $B(t,s)$ and $(H_x(\Delta), \|\cdot\|_x)$ be the Hilbert space with the reproducing kernel $B(t,s)$ on $\Delta \times \Delta$. If $0 < \widetilde{\theta}(a, \Delta) < \infty$, then

$$\left| \sqrt{\widetilde{\theta}(a+f, \Delta)} - \sqrt{\widetilde{\theta}(a, \Delta)} \right| \leq \|f\|_x / \sqrt{2|\Delta|}. \tag{19}$$

**Remark 4**. Lemma 3 is a version of Proposition 1.6 from the paper by Aurzada and Dereich (2011); relation (19) successfully supplements the original Lemma 1.

**Proposition 3.** The persistence exponent of process $S_H(t)$ has the following estimates

$$(1-H)H^{-1} \ln 1/(2H) \leq \widetilde{\theta}_{S_H} \leq 49 H^{-2}, \tag{20}$$

where the left inequality holds for $0 < H \leq e^{-2}/2$.

**Corollary**: *odd component of the fractional Brownian motion*.

Consider $w_H^-(t) = (w_H(t) - w_H(-t))/2$. Its dual stationary process $\widetilde{w}_H^-$ has the following correlation function:

$$\widetilde{B}_{w_H^-}(t) = (\cosh \frac{t}{2})^{2H} - (\sinh \frac{t}{2})^{2H}.$$

The exponent $\widetilde{\theta}_{w_H^-}$ exists because $\widetilde{B}_{w_H^-}(t) \geq 0$. By comparing $\widetilde{B}_{w_H^-}(t)$ with $\widetilde{B}_{w_H}(pt)$, Krug et al (1997) estimated the exponent as follows:

$$\widetilde{\theta}_{w_H^-} \geq \min((1-H)^2/H, \ (1-H)2^{1/(2H)-1}), \ 0 < H < 0.5, \tag{21}$$

$$\widetilde{\theta}_{w_H^-} \leq (1-H)^2/H, \quad 0.1549 < H < 0.5.$$

For small H these estimates are one-sided only. The following inequality

$$\widetilde{B}_{w_H^-}(2t) = (\cosh t)^{2H}(1 - (\tanh t)^{2H}) \geq (\cosh t)^{2H}(1 - |t|^{2H})_+ \geq \widetilde{B}_{S_H}(t)$$

and Proposition 3 immediately yield

$$\widetilde{\theta}_{w_H^-} \leq (7/H)^2/2, \quad 0 < H < 0.5.$$

**Remark 5.** A considerable difference in the behavior of $\widetilde{\theta}_{w_H^-}$ and $\widetilde{\theta}_{w_H} = 1 - H$ for small $H$ is expected. Heuristically this can be explained as follows. As $H \to 0$, the discrete processes $\widetilde{w}_H^-(k\Delta)$ and $\widetilde{w}_H(k\Delta)$ have different weak limits: $\{\xi_k\}$ and $\{(\xi_k - \eta)/\sqrt{2}\}$, respectively, where $\{\xi_k\}$ and $\eta$ are independent standard Gaussian variables. The probability (4) for the limiting processes are quite different:

$$P\{\xi_k < 0, k = 1 \div N\} = 2^{-N} \text{ and } P\{\xi_k - \eta \leq 0, k = 1 \div N\} = (N+1)^{-1}.$$

Unfortunately, this argument is insufficient to predict the behavior of $\widetilde{\theta}_{S_H}$ for small $H$, because the step $\Delta$ cannot be arbitrary and is a function of $H$.

**2.4 Khanin's problem.**

The survival exponent for fractional Brownian motion in the intervals $\Delta_T = (-T, T)$ is independent of the parameter $H$: $\theta_{w_H} = 1$. This interesting fact follows from both self-similarity of $w_H$ and the stationarity of its increments (Molchan, 1999).

In the case $H < 0.5$, variables $w_H(t)$ and $w_H(-t)$ are positive correlated. Therefore, a possible power-law asymptotics

$$P(w_H(t) < 1, -w_H(-t) < 1, t \in (0, T)) = T^{-\theta + o(1)}, \tag{22}$$

where we change sign before $w_H(t)$ for negative $t$ only, may have a radically different exponent compared with $\theta_{w_H} = 1$. The question of finding bounds on the exponent $\theta_{\chi_H}$ for the process

$$\chi_H(t) = sign(t)w_H(t), \quad \Delta_T = (-T,T)$$

was asked by K.Khanin. The next proposition contains a partial answer to this question.

**Proposition 4. 1.** In the case $0.5 \leq H < 1$, the exponent $\theta_{\chi_H}$ for $\Delta_T = (-T,T)$ exists and admits of the following estimates

$$1 < \theta_{\chi_H}(1-H)^{-1} \leq 2, \quad 0.5 \leq H < 1,$$

in addition $\theta_{\chi_{1/2}} = 1$.

2. Let $\underline{\theta}_{\chi_H}$ be the lower exponent in (22); then

$$\underline{\theta}_{\chi_H}(1-H)^{-1} \geq (H^{-1}-1) \wedge 2^{1/(2H)-1}, \quad 0 < H < 0.25,$$

$$\underline{\theta}_{\chi_H}(1-H)^{-1} \geq 2, \quad 0.25 < H \leq 0.5.$$

**Remark 6.** To clarify why $\theta_{\chi_H} / \theta_{w_H}$ is unbounded for small $H$ in the case $\Delta_T = (-T,T)$, we consider again the limiting sequence for $w_H(k\Delta)$ as $H \to 0$. This is $\{(\xi_k - \xi_0)/\sqrt{2}\}$, where the $\{\xi_k\}$ are independent standard Gaussian variables. The probability (1) for the limit sequence is

$$P\{\xi_k < \xi_0 + \sqrt{2}, |k| \leq N\} = (2N+1)^{-1}l(N),$$

where $l(N)$ is a slowly varying function, whereas for the limit sequence of $\chi_H(k\Delta)$ we have

$$P\{\xi_{-k} - \sqrt{2} < \xi_0 < \xi_k + \sqrt{2}, 0 < k \leq N\} \approx \sqrt{\pi}eN^{-1/2}\Phi(\sqrt{2})^{2N}, \quad (23)$$

where $\Phi(x)$ is the Gaussian distribution function. As in Remark 5, we have non-trivial exponential asymptotics where the threshold for $\{\xi_k\}$ is constant or bounded. Indeed, the event in (23) yields the inequality

$$|\xi_0| < \sqrt{2} + \max\left(\left|\sum_1^N \xi_{-k}\right|, \left|\sum_1^N \xi_k\right|\right)/N = \sqrt{2} + O_p(1)/\sqrt{N}.$$

**2.5 Explicit value of $\theta_x$.**

We have two explicit but isolated results for the fractional Brownian motion: $\theta_{w_H} = (1-H)$ for $\Delta_T = (0,T)$ and $\theta_{w_H} = 1$ for $\Delta_T = (-T,T)$. These results can be combined as follows:

**Proposition 5.** If $\Delta_T = (-T^\alpha, T)$, $0 \leq \alpha \leq 1$, then $\theta_{w_H} = \alpha H + (1-H)$.

**Remark 7.** The result is based on the following properties of the position $t_\Delta^*$ of maximum of $w_H(t)$ in $\Delta = [0,1]$: $t_\Delta^*$ has continuous probability density $f_\Delta^*(t)$ in $(0, 1)$ and $f_\Delta^*(t) \approx O(t^{-H})$ as $t \to 0$. In the case of multidimensional time, the behavior of $f_\Delta^*(t), \Delta = (0,1)^d$, near $t = 0$ is a key to the survival exponent for $w_H(t)$ in $\Delta_T = (-T^\alpha, T)^d$, if $0 < \alpha < 1$ and $H < 1$. By (2), $\theta_{w_H} = d$ in the case $\alpha = 1$, and $\theta_{w_H} = \alpha d$ in the degenerate case: $H = 1$.

## 3. Proofs.

### Proof of Proposition 3.

*Lower bound.* Let $\widetilde{w}_H(t)$ be the dual fractional Brownian motion with the parameter $H$, i.e., a Gaussian stationary process with correlation function $\widetilde{B}_{w_H}(t) = \cosh(Ht) - 0.5(2\sinh(t/2))^{2H}$. We prove in the Appendix that for $0 < H \leq e^{-2}/2$,

$$\widetilde{B}_{w_H}(pt) \geq B_{S_H}(t), \quad p = -H^{-1}\ln(2H). \tag{24}$$

Applying Slepian's lemma, one has $\widetilde{\theta}_{S_H} \geq p(1-H)$ because $\widetilde{\theta}_{w_H} = (1-H)$.

*Upper bound.* The random variable $\eta = \int_0^1 S_H(t)dt$ corresponds to an element $f_\eta(t)$ of the Hilbert space, $H_S(\Delta), \Delta = (0,1)$, with the reproducing kernel $B(t,s) = 1 - |t-s|^{2H}$. By definition of $H_x(\Delta)$, we have

$$f_\eta(t) = ES_H(t)\eta = 1 - (t^{1+2H} + (1-t)^{1+2H})/(1+2H),$$

$$\|f_\eta\|_S^2 = E\eta^2 = H(3+2H)(1+H)^{-1}(1+2H)^{-1}.$$

It is easy to see that $f_\eta(0) \leq f_\eta(t) \leq f_\eta(1/2)$. Therefore,

$$H < f_\eta(t) < 2H\ln(2e) \quad \text{and} \quad 4H/3 < \|f_\eta\|_S^2 < 3H. \tag{25}$$

Let $m_H$ be the median of the random variable $M = \max\{S_H(t), t \in \Delta)$, where $\Delta = (0,1)$. Then

$$0.5 = P(S_H(t) < m_H, t \in \Delta) < P(S_H(t) < m_H H^{-1} f_\eta(t), t \in \Delta),$$

because $H^{-1}f_\eta(t) > 1$. Setting $x(t) = S_H(t)$ in Lemma 2 and using notation (17), one has

$$\widetilde{\theta}(m_H H^{-1} f_\eta(t), \Delta) < \ln 2$$

and

$$\sqrt{\widetilde{\theta}(0,\Delta)} < \sqrt{\ln 2} + m_H H^{-1} \|f_\eta\|_S / \sqrt{2}.$$

Using Lemma 1 and the inequality $\|f_\eta\|_S < \sqrt{3H}$, we have

$$\widetilde{\theta}_{S_H} < \widetilde{\theta}(0,\Delta) < (\sqrt{\ln 2} + m_H \sqrt{1.5/H})^2.$$

It is well known (see, e.g., Lifshits, 1995) that $m_H < 4\sqrt{2}D(\Delta, \sigma/2)$, where $\sigma^2 = \max_\Delta ES_H(t)$ and D is the Dudley entropy integral related to the semi-metrics on $\Delta$: $\rho^2(t,s) = E(S_H(t) - S_H(s))^2$.

In our case $\rho(t,s) = \sqrt{2}|t-s|^H$, $\sigma = 1$ and therefore

$$m_H < c_H/\sqrt{H},$$

where

$$c_H = 4\sqrt{(1-H)\ln 2} + 2^{3-H}\sqrt{\pi}\Phi(-\sqrt{1-H}\ln 4) < 5.36 \quad, \qquad H < 1/2,$$

and $\Phi(x)$ is the standard Gaussian distribution. Hence,

$$\widetilde{\theta}_{S_H} < (\sqrt{\ln 2} + 5.36\sqrt{1.5}/H)^2 < (7/H)^2.$$

*Proof of Proposition 4.*

*Proposition 4.1.* In the case of $H \geq 0.5$, the process $\chi_H(t) = sign(t)w_H(t)$ has non-negative correlations on $R^1$. In the standard manner, this implies the existence of $\theta_{\chi_H}$ for $\Delta_T = (-T, T)$. More precisely, starting from a self-similar 2-D process $x(t) = (w_H(t), -w_H(-t))$ on $R^1_+$, we consider the dual 2-D stationary process $\widetilde{x}(t) = x(t)\exp(-Ht)$ whose correlation matrix has positive elements. Therefore, by Li and Shao (2004), we conclude that the exponent $\widetilde{\theta}_{\chi_H}$ for $\widetilde{x}(t)$ exists.

*Equality $\widetilde{\theta}_{\chi_H} = \theta_{\chi_H}$ for $\Delta_T = (-T, T)$.* We will use Lemma 1. By the relation $\chi_H(t) = sign(t)w_H(t)$, the map $\varphi(t) \mapsto sign(t)\varphi(t)$ is an isometry between the reproducing kernel Hilbert spaces $H_{\chi_H}(\Delta_T)$ and $H_{w_H}(\Delta_T)$ associated with $\chi_H(t)$ and $w_H(t)$ on $\Delta_T = (-T, T)$, respectively. To prove that the dual exponents are equal, it is enough to find $\varphi(t) \in (H_{w_H}(R^1), \|\cdot\|_{w_H})$ such that $\text{sgn}(t)\varphi(t) \geq 1$ for $|t| \geq 1$. We can use

$$\varphi(t) = \text{sgn}(t)\min(|t|, 1) = \int (e^{it\lambda} - 1)\frac{\sin\lambda}{i\pi\lambda^2}d\lambda,$$

because
$$\|\varphi\|^2_{w_H} = k_H \int \frac{(\sin\lambda)^2}{(\pi\lambda^2)^2}|\lambda|^{1+2H}d\lambda < \infty,$$

(see Molchan and Khokhlov, 2004).

*Estimation of $\theta_{\chi_H}$, $H > 1/2$.* Since $E\chi_H(t)\chi_H(s)) \geq 0$ for any $t, s$, we have, by Slepian's lemma,

$$p_T := P(w_H(t) < 1, -w_H(-t) < 1, t \in (0,T)) \geq [P(w_H(t) < 1, t \in (0,T))]^2.$$

Using (2), one has $\theta_{\chi_H} \leq 2(1-H)$.

Obviously, $p_T \leq P(w_H(t) < 1, t \in (0,T))$. Therefore, $\theta_{\chi_H} \geq (1-H)$ for any H.

*Proposition 4.2.* Let $0 < H \leq 1/2$, then $Ew_H(t)(-w_H(-s)) \leq 0$ for $t, s > 0$. Hence,

$$p_T \leq [P(w_H(t) < 1, t \in (0,T))]^2 \text{ and } \theta_{\chi_H} \geq 2(1-H).$$

Finally,

$$p_T \leq P(w_H(t) - w_H(-t) < 2, t \in (0,T)) = P(w_H^-(t) < 1, t \in (0,T)).$$

But then, $\theta_{\chi_H} \geq \theta_{w_H^-}$ for all $H$. If $\theta_{w_H^-} = \widetilde{\theta}_{w_H^-}$, then we get the lower bound of $\widetilde{\theta}_{\chi_H}$ for $0 < H \leq 1/4$.

*The equality* $\theta_{w_{\bar{H}}} = \tilde{\theta}_{w_{\bar{H}}}$. Let $H_{w_{\bar{H}}}(\Delta)$ and $H_{w_H}(\Delta)$ be the reproducing kernel Hilbert spaces associated with $w_{\bar{H}}(t)$ and $w_H(t)$, respectively. By the definition of $w_{\bar{H}}(t)$, the map $(\varphi(t), t > 0) \mapsto (sign(t)\varphi(|t|), |t| < \infty)$ is an isometric embedding of $H_{w_{\bar{H}}}(R_+^1)$ in $H_{w_H}(R^1)$. To prove that the exponents are equal, it is enough to find $\varphi(t), t \geq 0$ such that $sign(t)\varphi(|t|) \in (H_{w_H}(R^1), \|\cdot\|_{w_H})$, $\varphi(t) \geq 1$ for $t \geq 1$, and $\|\varphi\|_{w_H} < \infty$. As we showed above, this can be $\varphi(t) = \min(t,1), t > 0$.

***Proof of Proposition 5.*** Consider the fractional Brownian motion in $\Delta_T = (-T^\alpha, T)$, $0 \leq \alpha \leq 1$. By lemma 1, we can focus on the exponent related to the position of the maximum of $w_H(t)$ in $\Delta_T$, $t^*_{\Delta_T}$.

*Distribution of* $t^*_\Delta$. We remind the main properties of the distribution function, $F^*(x)$, of $t^*_\Delta$ related to the normalized interval $\Delta = (0,1)$ (see Molchan,1999; Molchan and Khokhlov,2004):

- $F^*(x)$ has continuous density $f_\Delta^*(x) > 0, 0 < x < 1$ such that $(1-x)f_\Delta^*(x)$ decreases and $xf_\Delta^*(x)$ increases on $\Delta$;

- $F^*(x)$ have the following estimates:

$$x^{1-H}l^{-1}(x) \leq F^*(x) \leq x^{1-H}l(x), \tag{26}$$

where $l(x) = \exp(c\sqrt{-\ln x}), c > 0$.

Due to monotonicity of $(1-x)f_\Delta^*(x)$ and $xf_\Delta^*(x)$, one has

$$(1-x)f_\Delta^*(x) \leq x^{-1}\int_0^x (1-u)f_\Delta^*(u)du \leq x^{-1}F^*(x), \tag{27}$$

$$xf_\Delta^*(x) \geq x^{-1}\int_{xq}^x uf_\Delta^*(u)du \geq q(F^*(x) - F^*(xq)), 0 < q < 1. \tag{28}$$

By (26, 27),

$$f_\Delta^*(x) \leq x^{-H}l(x)(1-x)^{-1}. \tag{29}$$

Using (26, 28), one has

$$f_\Delta^*(x) \geq qx^{-H}l^{-1}(x)(1 - l(x)l(xq)q^{1-H}).$$

If we set $q^{1-H} = l^{-2}(x)/2$, then

$$f_\Delta^*(x) \geq qx^{-H}l^{-1}(x)/2 = c_H x^{-H} l^{-\nu_H}(x), \tag{30}$$

where $\nu_H = (3-H)/(1-H)$, $c_H = 2^{-(2-H)/(1-H)}$.

*Distribution of* $t^*_{\Delta_T}$. Let $T_1 = T_- + T$, where $T_- = T^\alpha$, then the following processes $w_H(T_1\tau - T_-) - w_H(-T_-)$ and $w_H(\tau)T_1^H$ on $\Delta = (0,1)$ are equal in distribution. Hence, $t^*_{\Delta_T}$ and $T_1 t^*_\Delta - T_-$ have the same distribution as well. Therefore,

$$p_T := P(|t^*_{\Delta_T}| \leq 0.5) = P(|t^*_\Delta - T_-/T_1| \leq 0.5/T_1) = T_1^{-1} f_\Delta^*((T_- + \varepsilon)/T_1), \quad (31)$$

where $|\varepsilon| \leq 0.5$. We have used here the existence and continuity of $f_\Delta^*(x)$.

*Exponent* $\breve{\theta}_{w_H}$. Suppose $\alpha = 1$. Then (31) implies $\lim_{T \to \infty} T p_T = 0.5 f_\Delta^*(0.5)$.

Let $\alpha < 1$, then $(T_- + \varepsilon)/T_1 = o(1)$ as $T \to \infty$, and (30, 31) give a lower bound of $p_T$:

$$T_1 p_T \geq c_H (a_T^+)^{-H} l^{-\nu_H}(a_T^+).$$

Here and below $a_T^\pm = (T_- \pm 0.5)/T_1$.

Using (29, 31), we get an upper bound on $p_T$:

$$T_1 p_T = f_\Delta^*((T^\alpha + \varepsilon)/T_1) \leq (a_T^-)^{-H} l(a_T^-) T_1/(T+1) \leq 2(a_T^-)^{-H} l(a_T^-).$$

By substituting $T_- = T^\alpha$, we have

$$\ln a_T^\pm = -(1-\alpha)\ln T + O(T^{-\beta}), \beta = \alpha \wedge (1-\alpha) \quad \text{and} \quad \ln l(a_T^\pm) = O(\sqrt{\ln T}).$$

Hence,

$$-\ln p_T = (1 - (1-\alpha)H))\ln T + O(\sqrt{\ln T}),$$

i.e., $\breve{\theta}_{w_H} = \alpha H + (1-H)$.

*The equality* $\breve{\theta}_{w_H} = \theta_{w_H}$. Consider the Hilbert space $(\mathrm{H}_{w_H}(R^1), \|\cdot\|)$ related to FBM and a function

$$\varphi(t) = \min(|t|,1) = \int (e^{it\lambda} - 1)\left(\frac{\sin \lambda/2}{\sqrt{2\pi}\lambda/2}\right)^2 d\lambda. \quad (32)$$

The standard spectral representation of the kernel $E w_H(t) w_H(s)$ and the representation (32) yield

$$\|\varphi\|^2 = k_H \int \left(\frac{\sin \lambda/2}{\sqrt{2\pi}\lambda/2}\right)^4 |\lambda|^{1+2H} d\lambda < \infty,$$

where $k_H = \int |e^{i\lambda} - 1|^2 |\lambda|^{-1-2H}$. Setting $\varphi_T := \{\varphi(t), t \in \Delta_T\}$, the desired statement follows from Lemma 1 because $\varphi_T \in (\mathrm{H}_{w_H}(\Delta_T), \|\cdot\|_T)$ and $\|\varphi_T\|_T \leq \|\varphi\|$.

**Appendix.**

*Relation (7):* $\widetilde{B}_{I_H}(t) \leq \widetilde{B}_{I_{1-H}}(t)$.

By (6), one has for small and large $t$

$$\widetilde{B}_{I_H}(t) = 1 - (1-H^2)t^2/2 + (2+4H)^{-1}t^{2+2H}(1+o(1)), \quad t \to 0, \tag{A1}$$

$$\widetilde{B}_{I_H}(t) = (1+H)(1+2H)^{-1}e^{-Ht}(1-e^{-t}) + 0.5(1+H)e^{-\overline{H}t}(1+O(e^{-t})), \quad t \to \infty, \tag{A2}$$

where $\overline{H} = 1-H$. Therefore, we have the following asymptotics for $\Delta(t) = \widetilde{B}_{I_H}(t) - \widetilde{B}_{I_{\overline{H}}}(t)$:

$$\Delta(t) = -(1-2H)t^2/2 + O(t^{2+2H}), \quad t \to 0,$$

$$\Delta(t) = -(1-2H)H(2+4H)^{-1}e^{-Ht} - (1-2\overline{H})\overline{H}(2+4\overline{H})^{-1}e^{-\overline{H}t} + O(e^{-t}), \quad t \to \infty.$$

These relations support (7) both for small and large enough t. To verify (7) in the general case, we consider the following test function: $(2+4H)(2+4\overline{H})\Delta(t)\exp(-1.5t)$. Using new variables: $x = \exp(-t)$, $\alpha = 1-2H$, the test function is transformed in a function $\psi$ on square $S = (0,1) \times (0,1)$. Namely, $\psi = U(x,\alpha) - U(x,-\alpha)$, where

$$U(x,\alpha) = (4-\alpha^2)x^{3/2}(3-\alpha)\int_0^x [(x-u)((1-u)^{1-\alpha} - u^{1-\alpha}) + u^{1-\alpha}]du.$$

We have to show that $\psi \leq 0$. It is easy to see that $\psi = 0$ at the boundary of S. By (A1, A2), $\psi \leq 0$ in a neighborhood of two sides of S: x=0 and x=1. The same is true for the other sides: $\alpha = 0$ and $\alpha = 1$ because

$$\frac{\partial \psi}{\partial \alpha}(x,0) = -4(1-x^2)\int_{1-x}^1 \ln(1/u)du < 0,$$

and

$$\frac{\partial \psi}{\partial \alpha}(x,1) = (1-x)x^{-1/2}f(x) > 0.$$

Here $\quad f(x) = -x(1-x) + x^3 \ln 1/x + (1-x^3)\ln 1/(1-x)$.

To verify $f(x) > 0, 0 < x < 1$, note that $f'(x) = 3x^2(1 + v + \ln v)$, where $v = (1-x)/x$. Obviously, $f'$ has a single zero in (0,1), i.e. $f$ has a unique extreme point. But $f(0) = 0 = f(1)$ and $f(x) > 0$ for small $x$. Therefore $f(x) \geq 0, 0 < x < 1$.

Numerical testing supports the desired inequality $\psi < 0$ for interior points of $S$.

*Comment.* Our preliminary numerical test was concerned with points of grid with step 0.005. The first derivatives of $\psi$ are uniformly bounded from above on $S$. This fact helps to find a final grid step to prove $\psi < 0$ for all interior points of $S$. The corresponding analysis is unwieldy and so is omitted.

**Relation (9):** $\widetilde{B}_{I_H}(t) \leq \widetilde{B}_{I_{1/2}}(pt)$, $H \geq 1/2$, $p = 2(1-H)$.

To verify the inequality $\Delta(t) = \widetilde{B}_{I_H}(t) - \widetilde{B}_{I_{1/2}}(2(1-H)t) \leq 0$, we consider the following test function:
$(2+4H)\Delta(t)\exp(-(1+H)t)$. Using (6, 8) and new variables $(x = \exp(-t), \alpha = 2H - 1) \in S = (0,1) \times (0,1)$, we will have the following representation for the test function:

$$\psi(x,\alpha) = (3+\alpha)(x + x^{\alpha+2}) - 1 - x^{\alpha+3} + (1-x)^{\alpha+3} - 3(\alpha+2)x^2 + (\alpha+2)x^{3-\alpha} \quad (A3)$$

One has $\psi(x,\alpha) \leq 0$ in a neighborhood f two sides of S: x=0 and x=1, because

$$\psi(x,\alpha) = -(\alpha+2)(3-\alpha)x^2/2 + O(x^{(\alpha+2)\wedge(3-\alpha)}) < 0, x \to 0,$$

$$\psi(x,\alpha) = -2\alpha(1-\alpha)(3-\alpha)(1-x)^2/2 + O((1-x)^3) \leq 0, x \to 1.$$

The same is true for other sides: $\alpha = 0$ and $\alpha = 1$.

*Side* $\alpha = 0$. One has $\psi(x,0) = 0$ and

$$\frac{\partial \psi}{\partial \alpha}(x,0) = (1-x)[x(1-x) + 3x^2 \ln x + (1-x)^2 \ln(1-x)] := (1-x)\varphi_3(x) \leq 0$$

because

$$\varphi_a(x) = x(1-x) + ax^2 \ln x + (1-x)^2 \ln(1-x) \leq 0, \ a > 1 \quad (A4)$$

To prove (A4), note that $\varphi_a(0) = \varphi_a(1) = 0$ and $\varphi_a(x) = ax^2 \ln x + O(x^2) \leq 0$ as $x \to 0$. Hence, (A4) holds if $\varphi_a(x)$ has unique extremum in (0,1). By

$$\varphi_a^{(4)}(x) = -2ax^{-2} - 2(1-x)^{-2} \leq 0,$$

we conclude that $\qquad \varphi_a''(x) = (3a+1) + 2a \ln x + 2\ln(1-x)$

is a concave function with two zeroes in (0,1), because $\varphi_a''(1/2) > 0$ and $\varphi_a''(x) \to -\infty$ as $x \to 0$ or 1.

It means that

$$\varphi_a'(x) = (a-1)x + 2ax \ln x - 2(1-x)\ln(1-x)$$

has two extremums in (0,1) only. But $\varphi_a'(0) = 0$, $\varphi_a'(1) = a - 1 > 0$, and $\varphi_a'(x) \leq 0$ for small $x$ because $\varphi_a''(x) \to -\infty$ as $x \to 0$. Hence $\varphi_a'(x)$ has unique zero in (0,1) and $\varphi_a(x)$ has unique extremum.

So, we prove that $\psi(x,\alpha) \leq 0$ for small $\alpha$.

*Side* $\alpha = 1$. Here $\psi(x,1) = 0$ and

$$\frac{\partial \psi}{\partial \alpha}(x,1) = (1-x)(3-x)x^2 \ln(1/x) + (1-x)^2[x + (1-x)^2 \ln(1-x)] \geq 0,$$

because $\qquad [x + (1-x)^2 \ln(1-x)] \geq x + (1-x)\ln(1-x) = -\int_0^x \ln(1-u)du \geq 0.$

Hence, $\qquad \psi(x,\alpha) = \psi_\alpha'(x,1)(\alpha-1)(1+o(1-\alpha)) \leq 0, \ \alpha \to 1.$

As a result $\psi(x,\alpha) \leq 0$ near the boundary of $S = (0,1) \times (0,1)$. Numerical testing supports the desired inequality $\psi < 0$ for interior of $S$ (see more in the Comment from the Appendix section 'Relation 7').

**Relation (10):** $\widetilde{B}_{I_H}(t) \geq \widetilde{B}_{I_{1/2}}(pt)$, $H \leq 1/2$, $p = 2\sqrt{(1-H^2)/3}$.

Let $\psi = (2+4H)(\widetilde{B}_{I_H}(t) - \widetilde{B}_{I_{1/2}}(pt))e^{-(1+H)t}$. By change of variables: $x = \exp(-t)$ and $\alpha = 2H$, we get a test function

$$\psi(x,\alpha) = (2+\alpha)(x+x^{\alpha+1}) - 1 - x^{\alpha+2} + (1-x)^{\alpha+2} - 3(\alpha+1)x^{1+(\alpha+p)/2} + (\alpha+1)x^{1+(\alpha+3p)/2}$$

on $S=(0,1)\times(0,1)$ and the relation between $p$ and $\alpha$:

$$3(p/2)^2 + (\alpha/2)^2 = 1.$$

One has

$$\psi(x,\alpha) = (2+\alpha)x^{1+\alpha} - 3(1+\alpha)x^{1+(p+\alpha)/2} + O(x^2)) \geq 0, x \to 0,$$

$$\psi(x,\alpha) = (1-x)^{2+\alpha} + O((1-x)^3) \geq 0, x \to 1.$$

In addition, $\psi(x,0) = x(2 - 3x^{3^{-1/2}} + x^{3^{1/2}}) \geq 0$

Finally, $\psi(x,1) = 0$ and

$$\frac{\partial \psi}{\partial \alpha}(x,1) = \bar{x}(x\bar{x} + 2x^2 \ln x + \bar{x}^2 \ln \bar{x}) = \bar{x}\varphi_2(x),$$

where $\bar{x} = 1-x$. By (A4), $\varphi_2(x) \leq 0$.

Therefore $\psi(x,\alpha) \leq 0$ near the boundary of $S = (0,1) \times (0,1)$. The numerical testing supports this conclusion for interior of $S$ (see more in the Comment from the Appendix section 'Relation 7').

**Relations (13, 14).**

Consider $\Delta(t) = \widetilde{B}_{I_{1/2}}(t) - \widetilde{B}_L(pt)$, where $\widetilde{B}_L(t) = 1/\cosh(t/2)$ and $\widetilde{B}_{I_{1/2}}(t)$ is given in (8). By the change of variables $x = e^{-t/2}$, we transform the test function $2(1+e^{-pt})\Delta(t)$ in a function $\psi$ on $(0,1)$ such that

$$\psi(x) = (3x - x^3)(1+x^{2p}) - 4x^p.$$

Taking into account the asymptotics of $\psi$ near 0, we come to a necessary condition for $\psi$ to be negative, namely: $p \leq 1$. Let $p=1$, then $\psi = -(1-x^2)^2 x \leq 0$, i.e. $4\theta_L \leq 1$.

*Case* $p>1$. In this case $\psi \geq 0$ as $x \to 0$. An additional condition on $p > 1$ we can get from the relation $\psi \geq 0$ as $x \to 1$. One has $\psi = xQ(x)$, where

$$Q(x) = (3-x^2)(1+x^{2p}) - 4x^{p-1}.$$

By $Q(0) = 3, Q(1) = Q'(1) = 0$, we have $Q(x) = (1-x)^2 P(x)$ and $P(1) = 0.5Q''(1) = 2(p^2 - 3)$. Thus $Q(1) \geq 0$ if $p^2 \geq 3$.

*Case* $p = 2$. Here, $P(x)$ is a polynomial, $P(x) = 3 + 2x - 2x^3 - x^4$, and $P''(x) = -12x(1+x) \leq 0$, i.e. $P(x)$ is a concave function with $P(0) = 3, P(1) = 2$. Therefore $P(x) \geq 0$ and as a result, $4\theta_L \geq 1/p = 0.5$.

*Consider* $p = \sqrt{3}$. One has $Q(x) = 8(1-x)^3(1+o(1)), x \to 1$ and $Q(0) = 3 > 0$. Therefore $Q(x) \geq 0$, if $Q(x)$ is convex, i.e. $Q''(x) \geq 0$. To verify this property, note that
$$0.5x^2 Q''(x) = 2(3p-5)x^{p-1} + 3(6-p)x^{2p} - x^2 - (7+3p)x^{2+2p}$$
$$= (7+3p)x^{2p}(1-x^2) + \rho x^{p-1} + (1-\rho)x^{2p} - x^2 := \varphi(x),$$
where $\rho = 6p - 10$.

Obviously, $\varphi(x) \geq 0$ if $\rho x^{p-1} - x^2 \geq 0$. This holds for $0 < x < x_0 = 0.478$.

For $x > x_0$,
$$\rho x^{p-1} + (1-\rho)x^{2p} - x^2 \geq (\rho + (1-\rho)x_0^{p+1})x^{p-1} - x^2.$$

The right part here is positive for $x < 0.55$, i.e. $\varphi(x) \geq 0$ for $x \leq 0.5$.

Let $x > 0.5$. Then
$$\varphi(x) \geq (7+3p)2^{-2p}(1-x^2) + \rho x^{\alpha-1} + (1-\rho)x^{2p} - x^2$$
$$= C - (C+1)x^2 + \rho x^{p-1} + (1-\rho)x^{2p} := u(x),$$
where $C = (7+3p)2^{-2p}$. We have $u(0) = C, u(1) = 0$ and
$$u'(x) = -2(C+1)x + \rho(p-1)x^{p-2} + 2(1-\rho)px^{2p-1}$$
$$= -(C+1-2(1-\rho)px^{2p-2})x - ((C+1)x^{3-p} - \rho(p-1))x^{p-2}.$$
It is easy to see, that both terms in parentheses are positive on (0.5, 1).
Thus, $u(x)$ decreases to $u(1) = 0$. This means that $Q''(x) \geq 0$. Q.E.D.

**Relation (24):** $\widetilde{B}_{W_H}(pt) \geq B_{S_H}(t)$, $pH = -\ln(2H)$, $0 < H < e^{-2}/2$.

The difference of the correlation functions is the following
$$\Delta(t) = (cosh(Hpt) - 0.5(2\sinh(pt/2))^{2H}) - (1 - |t|^{2H})_+.$$
**Let** $t > 1$, then $\Delta(t) = \widetilde{B}_{W_H}(pt) \geq 0$.

**Let** $2H < t < 1$. It is enough to show that the first term, $\varphi$, in the following representation
$$\Delta(t) = [0.5e^{-Hpt} - 1 + t^{2H}] + 0.5e^{Hpt}(1 - (1 - e^{-pt})^{2H}) := \varphi + R$$

is non-negative. Setting $Hp = -\ln(2H)$, $\alpha = 2H$ one has
$$\varphi(t) = 0.5\alpha^t + t^\alpha - 1.$$

Let us show that $\varphi$ is decreasing. Then $\varphi$ is positive because $\varphi(1) = \alpha/2$.
We have
$$\varphi'(t) = \alpha^t(-05\ln(1/\alpha) + \psi(t)),$$
where $\psi(t) = \alpha^{1-t}/t^{1-\alpha}$. The function $\psi(t)$ has a single extreme point in the interval: $t^* = (1-\alpha)/\ln(1/\alpha)$. But $\psi(t^*) = \min$, because $\psi(t)$ decreases near $t = \alpha$:
$$\psi(\alpha) = 1 \quad \text{and} \quad \psi'(\alpha) = (\alpha\ln(e/\alpha) - 1)/\alpha \leq 0 \text{ for } 0 < \alpha < 1.$$

Hence, $\psi(t) \leq \max(\psi(\alpha), \psi(1)) = 1$. As a result,
$$\varphi'(t) \leq \alpha^t(-05\ln(1/\alpha) + 1) \leq 0.$$
The last inequality holds for $0 < \alpha < e^{-2}$.

So, we have,
$$\Delta(t) \geq 0, \ 2H < t < 1 \ \text{ for } 0 < \alpha < e^{-2}.$$

**Let** $0 < t < 2H$. Use

$$\Delta(t) = \cosh(Hpt) - 1 + t^{2H}[1 - 0.5(2t^{-1}\sinh(pt/2))^{2H}],$$

then $\Delta(t) \geq 0$ if
$$2^{1/(2H)} \geq \max_{(0,2H)}(2t^{-1}\sinh(pt/2)) = H^{-1}\sinh(pH) = (2H)^{-2} - 1.$$
This inequality holds for $0 < 2H < 1/4$.

Putting the above inequalities together yields (24) for $2H \leq e^{-2} \wedge 1/4$.

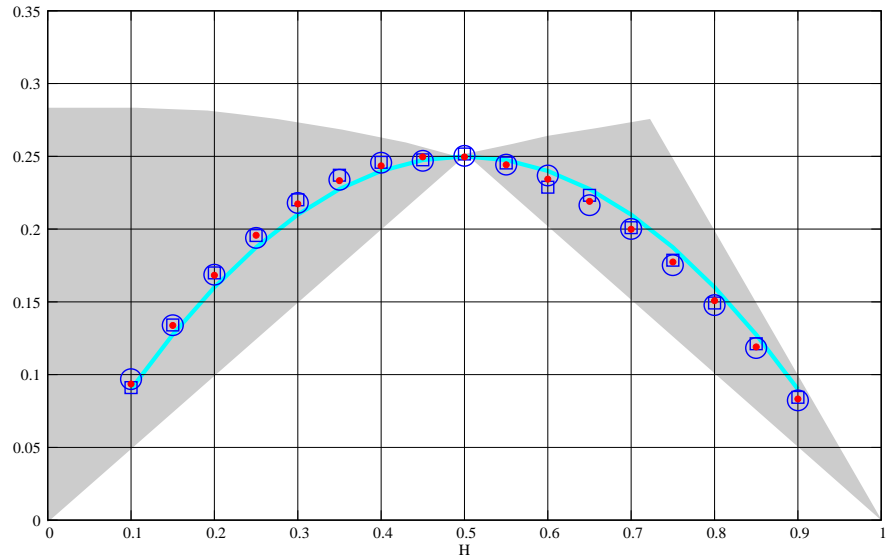

**Figure 1**

The survival exponents $\widetilde{\theta}_{I_H}$ for the integrated fractional Brownian motion in $\Delta_T = (-T, T)$: hypothetical values (*parabolic line*), empirical estimations (*small circles, squares*), and theoretical bounds (*shaded zone* given by Proposition **1(b,c)**).

The empirical exponents are based on the model (12, $\alpha(H) = 0$) in three time intervals of $\widetilde{T} = \ln T$:
$\ln(1/\varepsilon) \leq \widetilde{T}(1-H)H \leq \ln(10/\varepsilon)$ where $\varepsilon = 0.01$, 0.003, and 0.001 (see more in Molchan and Khohlov, 2003).